\documentclass[11pt]{amsart}

\usepackage{latexsym}
\usepackage{mathtools}
\usepackage{amssymb}
\usepackage{amsmath}
\usepackage{xcolor}
\usepackage[dvips]{graphicx}
\newlength{\cqfd}
\def\R{\mathbb{R}}

\newcommand{\T}{\Pi}

\usepackage[english,francais]{babel}
\newcommand{\into}{\int_{\Omega}}
\newcommand{\he}{h_\varepsilon}
\newcommand{\ue}{\bu_\varepsilon}

\newcommand{\eps}{\varepsilon}

\newcommand{\bu}{\bar{u}}
\newcommand{\udiv}{\, \mathrm{div}}

\newtheorem{theorem}{Theorem}[section]
\newtheorem{lemma}[theorem]{Lemma}
\newtheorem{e-proposition}[theorem]{Proposition}

\newtheorem{e-definition}[theorem]{Definition\rm}
\newtheorem{remark}{\it Remark\/}

\def\og{\leavevmode\raise.3ex\hbox{$\scriptscriptstyle\langle\!\langle$~}}
\def\fg{\leavevmode\raise.3ex\hbox{~$\!\scriptscriptstyle\,\rangle\!\rangle$}}

\begin{document}
% place in the next line the header (rubrique) chosen for your article,
% if you know it (you can also have 2, format : Header1/Header2
\centerline{}
%\begin{frontmatter}

% Title, authors and addresses

% use the thanksref command within \title, \author or \address for footnotes;
% use the ead command for the email address,
% and the form \ead[url] for the home page:
% \title{Title\thanksref{label1}}
% \thanks[label1]{}
% \author{Name\thanksref{label2}}
% \ead{email address}
% \ead[url]{home page}
% \thanks[label2]{}
% \address{Address\thanksref{label3}}
% \thanks[label3]{}
\selectlanguage{english}
\title[Quantitative estimates for Mean Field limits]{On Mean Field Limit  and Quantitative Estimates \\ with a Large Class of Singular Kernels: \\
        Application to the Patlak-Keller-Segel Model}

% use optional labels to link authors explicitly to addresses:
% \author[label1,label2]{}
% \address[label1]{}
% \address[label2]{}
% The [label1] can be suppressed if there is only one address for all authors

\selectlanguage{english}
\author{Didier Bresch}
\address{LAMA -- UMR5127 CNRS, Bat. Le Chablais, Campus Scientifique, 73376 Le Bourget du Lac, France}
\email{Didier.Bresch@univ-smb.fr}
\author{Pierre--Emmanuel Jabin}
\address{CSCAMM and Dept. of Mathematics, University of Maryland, 
    College Park, MD 20742, USA}
\email{pjabin@cscamm.umd.edu}
\author{Zhenfu Wang},
\address{Dept. of Mathematics, University of Pennsylvania, Philadelphia, PA 19104, USA} 
\email{zwang423@math.upenn.edu}
\thanks{D. Bresch is partially supported by SingFlows project, grant ANR-18-CE40-0027}
\thanks{P.--E. Jabin is partially supported by NSF DMS Grant 161453, 1908739 and NSF Grant RNMS (Ki-Net) 1107444}

%\address[authorlabel3]{Address3}

% If you know the dates of reception, and acceptation you can put them now;
%  idem the name of the person presenting the Note

%\medskip
%\begin{center}
%{\small Received *****; accepted after revision +++++\\
%Presented by }
%\end{center}
\maketitle
\begin{abstract}
\selectlanguage{english}
% Text of abstract in English
  In this note, we  propose a new relative entropy  combination of the methods developed by P.--E. Jabin and  Z.~Wang [Inventiones (2018)] and by S. Serfaty [Proc. Int. Cong. of Math, (2018) and references therein]  to treat more general kernels in mean field limit theory.  This new relative entropy may be understood as introducing appropriate weights in the relative entropy developed by P.-E. Jabin and Z. Wang (in the spirit of what has been recently developed by D.~Bresch and P.--E. Jabin [Annals of Maths (2018)])  to cancel the more singular terms involving the divergence of the flow. As an example, a full rigorous derivation (with quantitative estimates) of the  Patlak-Keller-Segel model  in some subcritical regimes is obtained. Our new relative entropy allows to treat singular potentials which combine large  smooth part, small attractive singular part and large repulsive singular part.
  %PN{\it To cite this article: A. Name1, A. Name2, C. R. Acad. Sci. Paris, Ser. I 340 (2005).}

\medskip
%\selectlanguage{francais}

% Text of abstract in French
\noindent{\bf R\'esum\'e:} {\it Limite champ moyen pour des noyaux plus g\'en\'eraux.}

\smallskip

\noindent  Dans cette note, on propose une nouvelle entropie relative combinant les m\'ethodes d\'evelopp\'ees par P.--E. Jabin et Z.~Wang [Inventiones (2018)] et par S. Serfaty [see review in Proc. Int. Con og Math (2018) et r\'ef\'erences] pour traiter des noyaux plus g\'en\'eraux en théorie de la limite champ moyen.
Cette nouvelle entropie relative consiste en l'introduction d'une famille de poids appropri\'es dans l'entropie relative d\'evelopp\'ee par  P.-E. Jabin et Z. Wang (dans le m\^eme esprit du travail r\'ecent de D.~Bresch et P.--E. Jabin [Annals of Maths (2018)]) pour compenser les termes les plus singuliers 
qui font intervenir la divergence du champ de vitesse. Comme exemple,  une preuve avec estimation quantitative de la limite champ moyen vers le mod\`ele de Patlak-Keller-Segel en r\'egime sous-critique est obtenue.   Notre m\'ethode permet de couvrir dses potentiels singuliers qui peuvent combiner une partie r\'eguliere, une petite partie singuli\`ere attractive et une grande partie singulière r\'epulsive.

%PN{\it Pour citer cet article~: A. Name1, A. Name2, C. R. Acad. Sci. Paris, Ser. I 340 (2005).}
\end{abstract}
%\end{frontmatter}

% now the Version française abrégée, if it exists
%\selectlanguage{francais}
\section*{Version fran\c{c}aise abr\'eg\'ee}
Dans cette note,  on pr\'esente une entropie relative permettant d'encoder quantitativement la limite champ moyen $N\to +\infty$ pour un syst\`eme de particules avec noyaux singuliers pouvant contenir
une partie r\'eguli\`ere qui peut \^etre grande, une partie singuli\`ere attractive sous critique, une grande partie singuli\`ere r\'epulsive.  Une preuve avec estimation quantitative de la limite champ moyen vers le mod\`ele de Patlak-Keller-Segel en r\'egime sous-critique est par exemple obtenu en consid\'erant le potentiel attractive de Poisson en dimension 2.
Par souci de simplicit\'e, on se restreint dans cette note au domaine p\'eriodique $\T^d$ avec $d=1,2,3$.
On consid\`ere la limite champ moyen $N\to +\infty$ pour le syst\`eme suivant \`a $N$ particules
\begin{equation}\label{sysfr}
d X_i = \frac{1}{N} \sum_{j\not = i} K(X_i-X_j)  dt+  \sqrt{2\sigma}  d B_i, \quad i=1, 2, \cdots, N, 
\end{equation}
o\`u les $B_i$ sont des mouvements Browniens ou des processus de Wiener ind\'ependants
et on consid\`ere un champ de vitesse type flot-gradient avec un noyau singulier donn\'e par
 \begin{equation}\label{kernelfr}
 K= -\nabla V.
 \end{equation}
Comme exemple en lien avec le mod\`ele de Patlak-Keller-Segel que l'on rencontre en biologie,
on donne le potentiel attractive de Poisson en dimension 2:
 \begin{equation}\label{logfr}
 V= \lambda \log |x| + V_e(x)
 \end{equation}
 o\`u $V_e$ une correction r\'eguli\`ere permettant \`a $V$ d'\^etre p\'eriodique. Le but est alors
 de donner une estimation quantitative de la convergence vers la solution du syst\`eme de
McKean-Vlasov
\begin{equation} \label{MFlimit0fr}
\partial_t \bar \rho + \udiv_x \left( \bar  \rho\,  K \star_x \bar \rho \right) = \sigma \Delta_x \bar \rho,  \quad \hbox{ with } \quad K = - \nabla V,\quad \bar\rho(t=0,x)=\bar\rho^0\in {\mathcal P}(\T^d)
\end{equation}
impliquant par voie de fait une preuve de propagation du chaos.  Il existe diff\'erents moyen
de comparer  \eqref{sysfr} avec la limite $\bar\rho$ donn\'ee par  \eqref{MFlimit0fr}. 
On suit ici \cite{JW1,JaWa} en introduisant la loi  $\rho_N(t,x_1,\ldots,x_N)$ satisfaisant l'\'equation
de Liouville 
  \begin{equation}
    \partial_t \rho_N
  + \sum_{i=1}^N {\rm div}_{x_i} \Bigl(\rho_N \frac{1}{N} \sum_{j \ne i }^N K(x_i-x_j)\Bigr) 
  = \sigma \sum_{i=1}^N  \Delta_{x_i} \rho_N,  \qquad  \rho_N\vert_{t=0} = \rho_N^0
   \hbox{ avec } \int_{\Pi^{dN}} \rho_N^0 = 1.\label{liouville0fr}
  \end{equation}

 \noindent \`A titre d'exemple, dans le cas du potentiel  de Poisson en dimension $d$, on obtient le r\'esultat suivant 

\medskip

\begin{theorem} Supposons $\rho_N\in L^\infty(0,T;\;L^1(\T^{Nd}))$ est solution entropique de \eqref{liouville0fr} avec la donn\'ee initiale  $\rho_N(t=0)=\bar\rho^{\otimes^N}(t=0)$, et le potentiel d'interaction $V$ donn\'e par \eqref{logfr}. Supposons que  $\bar\rho \in L^\infty(0,T;W^{2,\infty}(\Pi^d))$ est solution de \eqref{MFlimit0fr}  avec $\inf \bar\rho >0$. Supposons finalement que $\lambda<c_d\,\sigma$ pour une constante explicite $c_d$ dependant seulement de la dimension ($c_d=2$ si $d=2$ par exemple). Alors il existe une constante $C$ et un exposant $\theta>0$ ind\'ependant de $N$ tel que pour tout $k$  fix\'e:
\[
\|\rho_{N,k}-\bar\rho^{\otimes^k}\|_{L^\infty(0,T;\;L^1(\T^{kd}))}\leq C\,k^{1/2}\,N^{-\theta}
\]\label{PKSestimatefr}
o\`u $\rho_{N,k}$ est la marginale du syst\`eme au rang $k$ c'est-\`a-dire
$$\rho_{N,k}(t,x_1,\cdots,x_k) = \int_{\Pi^{(N-k)d}}
     \rho_N(t,x_1, \cdots,x_N) \, d x_{k+1}\cdots d x_{N}.
$$
\end{theorem}
Ce r\'esultat pr\'ecise \cite{FoJo} en donnant la premi\`ere justification rigoureuse de la propagation du chaos avec un taux de convergence explicite pour un potentiel singulier attractif.  Un r\'esultat plus g\'en\'eral peut  en fait \^etre obtenu permettant de consid\'erer des potentiels
$V$ du type
\begin{equation}
V(-x)=V(x),\qquad V= V_a+ V_r + V_s,\label{Vsumfr}
\end{equation}  
o\`u  $V_a$  un potentiel singulier attractif contraint, $V_r$  un potentiel singulier r\'epulsif assez g\'en\'eral et $V_s$ est un potentiel suffisamment r\'egulier. Le lecteur int\'eress\'e peut consulter Theorem \ref{Main} de la version anglaise pour plus de 
d\'etails concernant notamment les hypoth\`eses sur $V_s$, $V_r$ et $V_a$.

\smallskip

La preuve est bas\'ee sur une nouvelle entropie relative \eqref{modulateddef} qui est une combinaison des m\'ethodes d\'evelopp\'ees par P.--E. Jabin et Z.~Wang [Inventiones (2018)] et par S. Serfaty [see review in Proc. Int. Con og Math (2018) and references therein]. Cette nouvelle entropie relative peut aussi s'interpr\^eter comme une energie libre modifi\'ee du syst\`eme consistant en l'introduction d'une famille de poids appropri\'ee (\'equlibre de Gibbs $G_N$ et distribution correspondante $G_{\overline \rho_N}$
d\'efinis dans \eqref{Gibbs}) dans l'entropie relative d\'evelopp\'ee par  P.-E. Jabin et Z. Wang (dans le m\^eme esprit du travail r\'ecent de D.~Bresch et P.--E. Jabin [Annals of Maths (2018)]) pour compenser les termes les plus singuliers qui font intervenir la divergence du champ de vitesse.

\selectlanguage{english}
% main text

\section{Introduction.}
\label{sec2}
\noindent
      In this note, we consider the limit when $N\to +\infty$
for the following system of $N$ particles, 
\begin{equation}\label{sys}
d X_i = \frac{1}{N} \sum_{j\not = i} K(X_i-X_j)  dt+  \sqrt{2\sigma}  d B_i, \quad i=1, 2, \cdots, N, 
\end{equation}
where the $B_i$ are independent Brownian Motions or Wiener processes.
For simplicity, we limit ourselves in this note to the periodic domain $\T^d$ and we specifically consider gradient flows with 
 the interaction kernel given by
 \begin{equation}\label{kernel}
 K= -\nabla V.
 \end{equation}
 A guiding example in this note is the {\em attractive} Poisson potential in dimension $2$
 \begin{equation}
   V=\lambda\,\log |x|+V_e(x),\label{logpotential}
 \end{equation}
 where $V_e$ is a smooth correction so that $V$ is periodic. Logarithmic potentials still play a critical role if the dimension $d>2$ and for this reason we will still consider potentials like \eqref{logpotential} in any dimension, even if there is no connection with the Poisson equation anymore.
 
Our main goal is to provide precise quantitative estimates 
for the convergence of \eqref{sys} towards the limit McKean-Vlasov PDE 
\begin{equation} \label{MFlimit}
\partial_t \bar \rho + \udiv_x \left( \bar  \rho\,  K \star_x \bar \rho \right) = \sigma \Delta_x \bar \rho,  \quad \hbox{ with } \quad K = - \nabla V,\quad \bar\rho(t=0,x)=\bar\rho^0\in {\mathcal P}(\T^d). 
\end{equation}
In the case where $V$ is given by \eqref{logpotential} and $d=2$ then \eqref{MFlimit} is the famous Patlak-Keller-Segel model, which is one of the first models of chemotaxis for micro-organisms. The potential $-V \star \bar \rho $ can then be seen as the concentration of some chemical (one has typically $V\leq 0$ here): From \eqref{logpotential}, one has that $\Delta V-V=2\pi\,\lambda\, \delta_0$ so that the chemical is produced by the population. Moreover \eqref{MFlimit} implies that the population follows the direction of higher chemical concentrations (more negative values of $V$).

\smallskip

Although it offers only a rough modeling of the biological processes involved in chemotaxis, the Patlak-Keller-Segel model is a good example of a singular attractive dynamics (all micro-organisms try to concentrate on a point) competing with the spreading effect due to diffusion. Because $V\sim \log |x|$ is actually critical in dimension $2$, Eq. \eqref{MFlimit} may actually blow-up and form a Dirac mass in finite time. Still in dimension $2$, one may exactly characterize that such a blow-up occurs if and only if $\lambda>4\,\sigma$ (since we normalize $\bar\rho^0$ to have total mass $1$, this corresponds to the classical $8\,\pi\,\sigma$ critical mass if instead $V$ is normalized).  We refer for instance to \cite{BlDolPer,DoSe,DoPe} and the references therein. 

\smallskip

Because of this potential singular behavior, a full rigorous derivation  of the Patlak-Keller-Segel model from the stochastic \eqref{sys} has remained elusive, in spite of recent progress in \cite{CaPe,FoJo} or \cite{GQ,HaSc}. The results in \cite{FoJo} for example prove that any accumulation point as $N\to \infty$ of the random empirical measure associated to the system \eqref{sys} is a weak solution to \eqref{MFlimit} provided that one is in the so-called very subcritical regime with $\lambda<\sigma$. While this provides the mean field limit, at least in some weak sense, it does not imply propagation of chaos.
  Of course potentials like \eqref{logpotential} are  only one example of singular interactions between particles for which the mean field limit remain poorly understood, especially in the stochastic cases.  

\smallskip

There are several ways to quantitatively compare \eqref{sys} with the limit $\bar\rho$ given by \eqref{MFlimit}. We follow here \cite{JW1,JaWa} by using the joint law $\rho_N(t,x_1,\ldots,x_N)$ of the process $X_1,\ldots,X_N$ which solves the Liouville or Kolmogorov forward equation 
  \begin{equation}
    \partial_t \rho_N
  + \sum_{i=1}^N {\rm div}_{x_i} \Bigl(\rho_N \frac{1}{N} \sum_{j \ne i }^N K(x_i-x_j)\Bigr) 
  = \sigma \sum_{i=1}^N  \Delta_{x_i} \rho_N,  \qquad  \rho_N\vert_{t=0} = \rho_N^0.\label{liouvilleN}
  \end{equation}
In addition to quantitative convergence estimates, Eq. \eqref{liouvilleN} also offers a straightforward manner to understand solutions to system \eqref{sys}, which is actually non-trivial when dealing with singular potentials $V$. We mostly defer to the coming \cite{BrJaWa1} for a complete discussion of the notion of solutions that we require and which we simply call {\em entropy solutions} throughout this note: Those roughly corresponds to solutions in the sense of distribution to \eqref{liouvilleN} that also satisfy appropriate entropy and energy bounds. 
  
 \smallskip
 
The joint law $\rho_N$ is compared to the chaotic law
$\bar \rho_N := \bar \rho^{\otimes^N}=\Pi_{i=1}^N \bar\rho(t,x_i)$, built from the limit $\bar \rho$.
Of course $\bar\rho_N$ cannot be an exact solution to \eqref{liouvilleN} but  instead solves
\begin{equation}
\partial_t \bar \rho_N + \sum_{i=1}^N \udiv_{x_i}  \Bigl(\bar \rho_N  \, K\star_x \bar \rho(x_i) \,  \Bigr) 
  = \sigma\, \sum_{i=1}^N \Delta_{x_i} \bar \rho_N.   \label{tensorRhoN}
   \end{equation}
As probability densities, both $\rho_N$ and $\bar\rho_N$ are initially normalized by
\begin{equation} \label{Norm}
 \int_{\T^{dN}} \rho_N\vert_{t=0} = 1 =  \int_{\T^{d}} \bar  \rho\vert_{t=0},
\end{equation}
which is formally preserved by either \eqref{liouvilleN} or \eqref{MFlimit}/\eqref{tensorRhoN}.
  The method leads in particular to direct estimates between $\bar\rho^{\otimes k}$ and any  observable or marginal of the system at a fixed rank $k$,
\[
\rho_{N,k}(t,x_1,\ldots,x_k)=\int_{\T^{(N-k)\,d}} \rho_N(t,x_1,\ldots,x_N)\,dx_{k+1}\ldots x_{N}.
\]
A corollary of the analysis sketched in this note and fully developed in the coming \cite{BrJaWa1} is a rigorous derivation of the Patlak-Keller-Segel system in some subcritical regimes with as an example the following result
\begin{theorem} Assume that $\rho_N\in L^\infty(0,T;\;L^1(\T^{Nd}))$ is an entropy solution to Eq. \eqref{liouvilleN} normalized by \eqref{Norm}, with initial condition $\rho_N(t=0)=\bar\rho^{\otimes^N}(t=0)$, and for the potential $V$ given by \eqref{logpotential}. Assume that  $\bar\rho \in L^\infty(0,T;W^{2,\infty}(\Pi^d))$ solves Eq. \eqref{MFlimit}  with $\inf \bar\rho >0$. Assume finally that $\lambda<c_d\,\sigma$ for some explicit constant $c_d$ depending only on the dimension ($c_d=2$ if $d=2$ for example). Then there exists a constant $C>0$ and an exponent $\theta>0$ independent of $N$
  s.t. for any fixed $k$ 
\[
\|\rho_{N,k}-\bar\rho^{\otimes^k}\|_{L^\infty(0,T;\;L^1(\T^{kd}))}\leq C\,k^{1/2}\,N^{-\theta}.
\]\label{PKSestimate}
\end{theorem}
Theorem \ref{PKSestimate} follows directly from Theorem \ref{Main} stated later in this note and the classical Csisz\'ar-Kullback-Pinsker inequality. The exponent $\theta$ could be made fully explicit and actually depends only on $c_d\,\sigma-\lambda$. Ideally, we would hope to have convergence whenever we are in the subcritical regime with no blow-up, which would correspond to $\lambda<4\,\sigma$ in dimension $d=2$. Unfortunately, our analysis is currently not precise enough for this and Theorem \ref{PKSestimate} has instead the requirement $\lambda<2\,\sigma$ if $d=2$, losing a factor $2$ with respect to the optimal conjecture. The limitation may derive from our use of general large-deviation like estimates; we will instead develop specific estimates in the coming \cite{BrJaWa1} and expect to extend the assumption to $\lambda<4\,\sigma$. Theorem \ref{PKSestimate} is nevertheless the first rigorous derivation of the propagation of chaos with an  explicit rate  for such singular attractive potentials.

\smallskip

In the rest of this note, we first introduce the new relative entropy with weights related to the Gibbs equilibrium $G_N$ and the corresponding distribution $G_{\bar\rho_N}$ given by \eqref{Gibbs} in Section \ref{NRE}. We then sketch  the main steps in the estimates together with some non-optimal assumptions on the kernels in Section \ref{tools} and we conclude with the main theorem in Section \ref{general}.

%%%%%%%%%%%%%%%%%%%%%%%%%%%%%%%%%%%%%%%%%%%%%%%%%%%%%%%5
\section{A new relative entropy}
\label{NRE} 
%%%%%%%%%%%%%%%%%%%%%%%%%%%%%%%%%%%%%%%%%%%%%%%%%%%%%%%%
%
The main idea of the method is to base quantitative estimates on the free energy which is the natural physical notion for stochastic systems such as \eqref{sys}. Such a modulated free energy can actually be written as a relative entropy between the joint law $\rho_N$, the chaotic law $\bar\rho_N$ and the corresponding Gibbs equilibria, leading to 
\begin{equation}
\begin{split}
& E_{N} (\frac{\rho_N}{G_N}\,|\;\frac{\bar \rho_N}{G_{\bar \rho_N}})    =\frac{1}{N} %\Bigl[
  \int_{\Pi^{dN}} \rho_N (t,X^N) \log \Bigl(\frac{\rho_N(t,X^N) }{G_N(X^N)}\frac{G_{\bar  \rho_N}(t,X^N)}{\bar \rho_N(t,X^N)}\Bigr) dX^N,   \\ 
%& \hskip2cm -  \int_{\pi^{dN}} \rho_N (t,X^N) \, dX^N
 %  +\int_{\Pi^{dN}}  \frac{\bar \rho_N(t,X^N)} {G_{\bar \rho_N}(t,X^N)}  G_N(t,X^N\bigr)
  %      dX^N\Bigr]
\end{split}\label{modulateddef}
\end{equation}
where we denote by $G_N$ the Gibbs equilibrium of the system \eqref{sys}, and by $G_{\bar\rho_N}$ the corresponding distribution where the exact field is replaced by the mean field limit according to the law $\bar\rho$, leading to
\begin{equation}
\begin{split}
  & G_N(t,X^N) = \exp \bigg(- \frac{1}{2N\sigma } \sum_{i \ne j} V (x_i-x_j)\bigg),\\
  & G_{\bar\rho}(t,x)=\exp \bigg(-\frac{1}{\sigma}\,V\star \bar \rho(x) +\frac{1}{2\,\sigma}\,\int_{\Pi^{d} } V\star \bar \rho \, \bar\rho \bigg),  \\
& G_{\bar \rho_N} (t,X^N) = \exp \bigg(-  \frac{1}{\sigma} \sum_{i=1}^NV\star \bar \rho(x_i) 
     +  \frac{N}{2\sigma} \int_{\Pi^{d} } V\star \bar \rho \, \bar\rho \bigg).
\end{split}\label{Gibbs}
\end{equation}

\medskip

\noindent One may decompose $E_N$ as follows
 \[
E_{N} (\frac{\rho_N}{G_N}\,|\;\frac{\bar \rho_N}{G_{\bar \rho_N}}) 
   = {\mathcal H}_N(\rho_N\vert \bar \rho_N) 
    + {\mathcal K}_{N} (G_N\vert G_{\bar \rho_N}) 
    %+ {\mathcal L}_{\bar \rho_N/G_{\bar \rho_N}} (G_N \vert G_{\bar \rho_N})
    ,
    \]
where
\[
  {\mathcal H}_N(\rho_N\vert \bar \rho_N)
   = \frac{1}{N} \int_{\Pi^{dN}} \rho_N(t,X^N) \log\Bigl(\frac{\rho_N(t,X^N)}{\bar \rho_N(t,X^N)}\Bigr)
   \, d X^N
   \]
 is exactly the relative entropy introduced in \cite{JW1,JaWa}  and 
\[
  {\mathcal K}_{N} (G_N\vert G_{\bar \rho_N}) 
    = - \frac{1}{N}  \int_{\Pi^{dN}} \rho_N(t,X^N) \log \bigl(\frac{G_N(t,X^N)}{G_{\bar \rho_N}(t,X^N)} \bigr)
    \, d X^N
    \]
  which  is the expectation of the modulated  energy on which the method developed in  \cite{Du,Se,Se1} is based.  Our approach may hence be seen {\em as a combination of the two methods respectively introduced in } \cite{JaWa} {\em and in}  \cite{Se}: 
 
 \medskip

\begin{itemize}  
\item I) The article \cite{JaWa} compares $\rho_N$ with $\bar \rho_N$ through the rescaled relative entropy $\mathcal{H}_N(\rho_N\vert \bar \rho_N)$.
This leads to a quantitative version of propagation of chaos for systems \eqref{sys} under the assumption that the kernel $K$ belongs to $\dot W^{-1,\infty} $ with  ${\rm div} \,K\in \dot W^{-1,\infty} $.  This in particular implied the first quantitative derivation of the incompressible Navier-Stokes system in dimension $2$ from point vortices with diffusion (in that case ${\rm div}\,K = 0$) where only qualitative limits were previously available in \cite{FHM,Osada86,Osada}.

More precisely the Navier-Stokes case corresponds to taking $K=\mbox{curl}\, V$ with $V$ still given by \eqref{logpotential} and at first glance, it may hence appear that the derivation of the Patlak-Keller-Segel with the kernel $K$ given by $K=-\nabla V$ should follow similarly. However the requirement ${\rm div} K\in \dot W^{-1,\infty}$ breaks the comparison: Of course ${\rm div}\,K=0$ if $K=\mbox{curl}\, V$ while ${\rm div}\, K\in \dot W^{-1,\infty}$ forces $V$ to be log-Lipschitz if $K=-\nabla V$. In general and because of this condition on the divergence, the result in \cite{JaWa} {\em underperforms in the gradient flow setting that we consider here}.   
The interested reader is further referred to \cite{Sa} for an introduction on this paper in the Bourbaki seminar.

\smallskip

\item II) The strategy followed in \cite{Se,Se1} is implemented for deterministic system ($\sigma=0$ in \eqref{sys}) and consists in using the modulated potential energy. Written at the scale of the empirical measure $\mu_N=\frac{1}{N}\,\sum_{i=1}^N \delta(x-X_i(t))$  the modulated energy reads
\[
\frac{1}{2}\int_{\T^{2\,d}\cap  \{x \ne y\}} V(x-y)\,(\mu_N(dx)-\bar\rho(x)\,dx)\,(\mu_N(dy)-\bar\rho(y)\,dy),
\]
of which ${\mathcal K}_{N} (G_N\vert G_{\bar \rho_N})$ would be the expectation w.r.t. $\rho_N$.
 In the deterministic setting, the corresponding method provides quantitative rate of convergence for any Riesz potential $V(x)=c\,|x|^{-\alpha}$ with $c>0$ and $\alpha<d$. This not only applies to gradient flow settings but also to some Hamiltonian flows such as $K=c\,\mbox{curl}\,|x|^{-\alpha}$ in dimension $2$. Because one may use $\alpha>d-2$, the result goes beyond Coulombian interactions and the type of singularity that could be handled in \cite{Hauray} for example.

The approach in \cite{Se,Se1} is however limited to repulsive potentials (the modulated energy is otherwise not positive) and to exact Riesz potentials. Though it may be possible to relax the later constraint, even adding a smooth perturbation to $V$ is challenging. 
\end{itemize}
\medskip
The hope is hence that adding ${\mathcal K}_{N} (G_N\vert G_{\bar \rho_N})$ with $\mathcal{H}_N(\rho_N\vert \bar \rho_N)$ will allow to combine the control on potentials with very singular but specific structure in \cite{Se} with estimates in \cite{JaWa} that require less structure but cover less singular potentials. In some senses the proposed modulated free energy may also be understood as introducing appropriate weights in the relative entropy used in \cite{JaWa} (in the spirit of what
 has been recently developed by some of the authors in  \cite{BrJa} in an other framework) to cancel the more singular terms  involving the divergence of the flow namely involving ${\rm div} \, K=-\Delta\,V$.

\medskip

\noindent The first step is obviously to prove the following explicit expression for the time evolution of the modulated free energy~$E_N$
\begin{e-proposition}
    Assume that $V$ is an even function,   and that $\rho_N$ is an entropy solution to \eqref{liouvilleN} and $\bar \rho$ solves \eqref{MFlimit}. Then the modulated free energy defined by \eqref{modulateddef} satisfies that
\[
\begin{split}
&\frac{d}{dt} E_{N}\left(\frac{\rho_N}{ G_{N}}\,\vert\;\frac{\bar\rho_N}{G_{\bar\rho_N}}\right) \le  -\frac{\sigma}{N}\,\int_{\Pi^{d\,N}} d\rho_N\,\left|\nabla\log \frac{\rho_N}{\bar\rho_N}-\nabla\log \frac{ G_N}{G_{\bar\rho_N}}\right|^2\\
&\ -\frac{1}{2} \int_{\Pi^{dN}} \int_{\Pi^{2\,d}\cap \{x\neq y\}}  \nabla V(x-y) \cdot \left(\nabla\log \frac{\bar\rho}{G_{\bar \rho}}(x)
- \nabla\log \frac{\bar\rho}{G_{\bar \rho}}(y)\right)\, (d\mu_N - d\bar\rho)^{\otimes 2} d\rho_N,
\end{split} 
\]\label{modulatedweak}
where $\mu_N=\frac{1}{N}\,\sum_{i=1}^N \delta(x-x_i)$ is  the empirical measure.  
\end{e-proposition}

\medskip

\noindent The strategy for obtaining propagation of chaos from Prop. \eqref{modulatedweak}
then relies on proving the two following points:
\begin{itemize}
\item Show that $E_N\left(\frac{\rho_N}{ G_{N}}\,\vert\;\frac{\bar\rho_N}{G_{\bar\rho_N}}\right)$ effectively controls the distance between $\rho_N$ and $\bar\rho_N$. Note that  $E_N$ is not a priori a positive quantity since $\rho_N/G_N$ and $\bar\rho_N/G_{\bar\rho_N}$ do not have the same mass. In this note our goal will be to obtain estimates like
  \[
E_N\left(\frac{\rho_N}{ G_{N}}\,\vert\;\frac{\bar\rho_N}{G_{\bar\rho_N}}\right)\geq \frac{1}{C}\,\mathcal{H}_N\left( \rho_N \,\vert\; \bar\rho_N  \right)-\frac{C}{N^\theta},
\]
for some constant $C$ and exponent $\theta$ independent on $N$.
\item Control the right-hand side in Prop. \eqref{modulatedweak} by $E_N\left(\frac{\rho_N}{ G_{N}}\,\vert\;\frac{\bar\rho_N}{G_{\bar\rho_N}}\right)$ or individually by $\mathcal{H}_N\left(\rho_N \,\vert\; \bar\rho_N \right)$ or $\mathcal{K}_N \left(  G_{N} \,\vert\; G_{\bar\rho_N} \right)$ plus some vanishing in $N$ correction of the form $C/N^\theta$. 
  \end{itemize}
  
\medskip

While it may first appear that the second point is more complex, it is not necessarily so. In fact in the Patlak-Keller-Segel case with $V$ given by \eqref{logpotential}, the right-hand side in Prop. \ref{modulatedweak} may be directly controlled by the large deviation like estimate provided in \cite{JaWa} and it is the positivity of the modulated free energy $E_N$ that is the critical point.

\medskip

As matter of fact, even for {\em smooth, attractive} potentials, it may occur that $E_N$ is negative with $\mathcal{K}_N$ negative and dominating $\mathcal{H}_N$. This issue forced the limitation to repulsive potentials in \cite{Se} but fortunately for our present method, a straightforward solution in the stochastic case is simply to remove from $V$ any large smooth part that could create issues.
 This leads to the definition of possibly corrections to the modulated free energy: 
For a given $W$ smooth enough and not necessarily connected to the dynamics, we define
\[
\mathcal{K}^W_N(G_N^W\,\vert\;G^W_{\bar\rho_N})=\frac{1}{N}\int_{\T^{dN}}  \rho_N \log \frac{G^W_{\bar \rho_N}}{G^W_N}\,dX^N,
\]
with similarly to the potential $V$
\[
\begin{split}
& G^W_N(X^N) = \exp \bigg(- \frac{1}{2N\sigma } \sum_{i\neq j} W (x_i-x_j)\bigg),\\ & G^W_{\bar\rho}(x)=\exp \bigg(-\frac{1}{\sigma}\,W\star \bar \rho(x) +\frac{1}{2\,\sigma}\,\int_{\Pi^{d} } W\star \bar \rho \, \bar\rho\bigg),  \\
& G^W_{\bar \rho_N} (t,X^N) = \exp \bigg(-  \frac{1}{\sigma} \sum_{i=1}^N W\star \bar \rho(x_i) 
     +  \frac{N}{2\sigma} \int_{\Pi^{d} } W\star \bar \rho \, \bar\rho \bigg).
\end{split}
\]
Of course such terms introduce corrections as well in Prop. \ref{modulatedweak} and we need to calculate
\begin{lemma}
  For $W\in W^{2,\infty}$, even, one has that
  \[\begin{split}
&\frac{d}{dt}\mathcal{K}^W_N(G_N^W\,\vert\;G^W_{\bar\rho_N})= \frac{d}{dt}\frac{1}{N}\int_{\Pi^{dN}} \rho_N \log \frac{G^W_{\bar \rho_N}}{G^W_N}\,dX^N\\
  &\quad=\int \rho_N\,\int_{\Pi^{2d}} \Delta W(x-y)\,(d\mu_N-d\bar\rho)^{\otimes^2}\,dX^N\\
  &\qquad  -\frac{1}{\sigma} \int \rho_N\,\int_{\Pi^{2d}} \int_{\Pi^d} \nabla W(z-x)\cdot\nabla V(z-y)\,\bar\rho(z)\,dz \,(d\mu_N-d\bar\rho)^{\otimes^2}\\
  &\qquad -\frac{1}{2\,\sigma} \int \rho_N\,\int_{\Pi^{2d}} \nabla V(x-y)\,(\nabla W\star(\mu_N-\bar\rho)(x)-\nabla W\star(\mu_N-\bar\rho)(y))\, (d\mu_N-d\bar\rho)^{\otimes^2}\\
  &\qquad-\frac{1}{2\,\sigma}\,\int \rho_N\int_{\Pi^{2d}} \nabla W(x-y)\,(\nabla V\star \bar \rho(x)-\nabla V\star\bar\rho(y))\,(d\mu_N-d\bar\rho)^{\otimes^2}.
  \end{split}
 \]\label{smoothcorrection}
\end{lemma}
Since $W$ is not connected to the dynamics, Lemma \ref{smoothcorrection} does not have as nice a structure as Prop. \ref{modulatedweak}. In particular the right-hand side now involves $\Delta W$ and for this reason we will only consider  $\mathcal{K}_N^W$ for smooth~$W$.
%%%%%%%%%%%%%%%%%%%%%%%%%%%%%%%%%%%%%%%%%%%%%%%%%%%%%%%%%%%%%%%%%%%%%%%%%%      %%%%%%%%%%%%%%%%%%%%%%%%%%%%%%%%%%%%%%%%%%%%%%%%%%%%%%%%%%%%%%%%%%%%%%%%%% 
\section{The main tools\label{tools}} 
%%%%%%%%%%%%%%%%%%%%%%%%%%%%%%%%%%%%%%%%%%%%%%%%%%%%%%%%%%%%%%%%%%%%%%%%%%
As Prop. \ref{modulatedweak} and Lemma \ref{smoothcorrection} emphasize, the heart of the proof consists in controlling terms of the form
\begin{equation}
\int_{\T^{dN}} d\rho_N\int_{\Pi^{2d}} f(x,y)\,(d\mu_N-d\bar\rho)^{\otimes^2}\label{partitiongeneral}
\end{equation}
in terms of $\mathcal{H}_N$ or $\mathcal{K}_N$. In the first case, we will rely on large deviation like estimates, sketched in subsection \ref{sec:largedeviations} while in the second case one needs to use the repulsive nature of the potential in the spirit of \cite{Se}.
%%%%%%%%%%%%%%%%%%%%%%%%%%%%%%%%%%%%%%%%%%%%%%%%%%%%%%%%%%%%%%%%%%%%%%%%%%
\subsection{Explicit basic large deviation estimate \label{sec:largedeviations}}
%%%%%%%%%%%%%%%%%%%%%%%%%%%%%%%%%%%%%%%%%%%%%%%%%%%%%%%

\medskip

\noindent We first recall the classical estimate, which was Lemma 1 in \cite{JaWa} for similar purposes 

\smallskip
\begin{lemma}
  For any $\rho_N,\;\bar\rho_N\in \mathcal{P}(\T^{dN})$, any test function $\psi \in L^\infty(\Pi^{dN})$, one has that for any  $\alpha>0$, 
  \[
\int_{\T^{dN}} \psi(X^N)\,d\rho_N\leq \frac{1}{\alpha} \,  \frac{1}{N}\int d\rho_N\,\log \frac{\rho_N}{\bar\rho_N}+ \frac{1}{\alpha} \frac{1}{N}\log \int_{\T^{dN}} e^{ \alpha   N\,\psi(X^N)}\,d\bar\rho_N. 
  \]\label{duality}
  \end{lemma}
Lemma \ref{duality} directly connects bounds on quantities like \eqref{partitiongeneral} to the relative entropy $\mathcal{H}_N$ and estimates on quantities that can be seen as partition functions
\begin{equation}\displaystyle 
\int_{\T^{dN}} e^{\displaystyle 
 N\,\int_{\Pi^{2d}} f(x,y)\,(d\mu_N-d\bar\rho)^{\otimes^2}}\,\bar\rho^{\otimes^N}\,dX^N. \label{partition}
  \end{equation}
It is hence natural to try to use large deviation type of tools to bound \eqref{partition}. Note however that our goals here are different from classical large deviation approaches: We do not try to calculate the limit as $N\to \infty$ of \eqref{partition} but instead to obtain bounds that are uniform in $N$.
  As a first example of such result, we recall the estimate from \cite{JaWa} which reads in the present context
 
\smallskip

\begin{theorem} {\rm (Theorem 4 in \cite{JaWa})}.
There exists a constant $c_d$ depending only on the dimension such that for any 
for any $\bar\rho\in L^\infty(\T^d)\cap \mathcal{P}(\T^d)$ and for any $f\in L^\infty(\T^{2d})$ with $f(x,y)=f(y,x)$ and
\[
\gamma:=c_d\,\left(\sup_{p\geq 1} \frac{\|\sup_y |f(x,y)|\,\|_{L^p_x}}{p}\right)^2 +c_d \left( \,\|\bar\rho\|_{L^\infty}\,\|f\|_{L^\infty} \right)^2 <1,
\]
then 
\[
\sup_{N \geq 2} \, \int_{\T^{dN}} e^{N\,\int_{\Pi^{2d}} f(x,y)\,(d\mu_N-d\bar\rho)^{\otimes^2}}\,\bar\rho^{\otimes^N}\,dX^N\leq \frac{2}{1-\gamma} \leq C < \infty.
\]\label{largedeviationJaWa}
\end{theorem}

\smallskip

\begin{remark}
  Observe that $\int_{\Pi^{2d}} f(x,y)\,(d\mu_N-d\bar\rho)^{\otimes^2}=\frac{1}{N^2}\, \sum_{i,j} \psi(x_i,x_j)$ for some $\psi$ which was the formulation used in {\rm \cite{JaWa}}. Moreover one has that
  \[
  \psi(x,y)=f(x,y)-\int f(x,z)\,\bar\rho(z)\,dz-\int f(z,x)\,\bar\rho(z)\,dz+\int f(z,z')\,\bar\rho(z)\,\bar\rho(z')\,dz\,dz',
  \]
  so that $\int \psi(x,y)\,\bar\rho(y)\,dy=0$ and of course by symmetry so does
  $\int \psi(x,y)\,\bar\rho(x)\,dx=0$ which were the two cancellations required in {\rm \cite{JaWa}}.  
\end{remark}
While the proof of Theorem \ref{largedeviationJaWa} in \cite{JaWa} was combinatorics, quite recently a probabilistic proof of this theorem was given in \cite{LimLuNo} and the method then applied to a  chemical reaction-diffusion model.

\smallskip

Theorem \ref{largedeviationJaWa} already provides some interesting control for our problem and in particular it can immediately bound for the right-hand side of Prop. \ref{modulatedweak} when the potential has logarithmic singularities like \eqref{logpotential}. For example if
$|\nabla V(x) |\leq
 {C}/{|x|}$,  then combining Theorem \ref{largedeviationJaWa} with Lemma \ref{duality} yields that
\begin{equation}
\begin{split}
  &\int_{\T^{dN}}\int_{\T^{2d}\cap\{x\neq y\}} \nabla V(x-y)\cdot \left(\nabla \log \frac{\bar\rho}{G_{\bar\rho}}(x)-\nabla \log \frac{\bar\rho}{G_{\bar\rho}}(y)\right)\,(d\mu_N-d\bar\rho)^{\otimes^2}\,d\rho_N\\
  &\qquad\leq C\,\mathcal{H}_N(\rho_N\,\vert\;\bar\rho_N)+\frac{C}{N},
\end{split}\label{righthandsideKS}
  \end{equation}
  provided that $\nabla \log \bar\rho,\; \nabla V \star \bar \rho \in W^{1,\infty}$ (recall the definition of $G_{\bar \rho}$ ) and for some constant $C$ depending on the corresponding norms. 

\smallskip
  
However Theorem \ref{largedeviationJaWa} cannot provide lower bounds on the modulated free energy $E_N$ when $V$ has attractive singularities like \eqref{logpotential}. Instead one can go back to the classical large deviation approaches as developed for example in \cite{Arous,Var} and obtain explicit bounds out of them.  Given a possibly unbounded functional $F:\;\mathcal{P} (\Pi^d) \to \R\cup \{\pm\,\infty \}$, we recall that the large deviation functional associated to $F$ is
\begin{equation}
I(F)=\max_{\mu\in \mathcal{P}(\Pi^d)} F(\mu)-\int_{\Pi^d} \mu\,\log \frac{\mu}{\bar\rho}\,dx.\label{largedeviationfunct}
\end{equation}
It is well known that if $F$ is continuous then at the limit $N\to \infty$, at first order the partition function given by \eqref{partition} behaves like $e^{N\,I(F)}$. We essentially have to quantify such estimates to make them uniform in $N$ and applicable to possibly unbounded $F$ where the classical theory typically requires $F$ to be continuous on measures.

We hence introduce an additional assumption ensuring the coercivity with respect to the relative entropy; namely we require that for some smooth convolution kernel $L_\eps$, one has that for some $k>0$,
\begin{equation}
\int_{\Pi^{dN}} |F(\mu_N)-F(L_\eps\star \mu_N)|\,\rho_N\,dX^N\leq C\,\eps^k+\frac{\alpha}{N}\int_{\Pi^{dN}} \rho_N\,\log \frac{\rho_N}{\bar \rho_N }.\label{regularityF}
\end{equation}
This allows to follow classical approaches and obtain for example

\smallskip

\begin{theorem}
  Assume that $\log\bar\rho\in W^{1,\infty}$, then there exists a constant $C$ depending only on $d$, $L$,  s.t. for any $F$ satisfying \eqref{regularityF}, one has that
  \[\begin{split}
  &\frac{1}{N}\log \int_{\Pi^{dN}} e^{N\,F(L_\eps\star\mu_N)}\,\bar\rho_N\,dX^N\leq I(F)\\
  &\qquad+\frac{C}{N^{1/(d+1)}\,\eps^{d/(d+1)}}\,(\log N+|\log \eps|+\|\log\bar\rho\|_{L^\infty})+C\,\eps\,\|\log\bar \rho\|_{W^{1,\infty}}.
  \end{split}
  \]
  Consequently, one has for any $\rho_N\in { \mathcal{P}}(\Pi^{dN})$ that for some exponent $\theta>0$ depending on $d$ and $k$
  \[\begin{split}
  &  \int_{\Pi^{dN}} F(\mu_N)\,\rho_N\,dX^N\leq I(F)+(\alpha +1)\, \mathcal{H}_N(\rho_N\,\vert \;\bar\rho_N)
  +C\,N^{-\theta}\,\left(1+\|\log\bar\rho\|_{W^{1,\infty}} + \mathcal{H}_N(\rho_N\,\vert \;\bar\rho_N)\right).
  \end{split}
  \]\label{largedeviation}
\end{theorem}
Theorem \ref{largedeviation} is our basic tool to control the modulated energy $\mathcal{K}_N$ in terms of the relative entropy $\mathcal{H}_N$. In that case, we take as $F$
\begin{equation}
F(\mu)= - \gamma \,\int_{\T^{2d}\cap\{x\neq y\}} V(x-y)\,(d\mu-d\bar\rho)^{\otimes^2},\label{specialF}
\end{equation}
for some $\gamma >0$.
  Provided that the attractive part of $V$ is not too singular so that it can be controlled by the entropy through \eqref{regularityF} and if that $I(F)=0$,  then we have that
\[
\mathcal{K}_N(G_N\,\vert\;G_{\bar\rho_N})\geq -\frac{\alpha+1}{ 2 \gamma \sigma }\,\mathcal{H}_N(\rho_N\,\vert\;\bar\rho_N) -\frac{C}{N^\theta}, 
\]
which will give the desired bound provided that ${(\alpha +1)}/{(2 \gamma  \sigma)}<1$.  
%%%%%%%%%%%%%%%%%%%%%%%%%%%%%%%%%%%%%%%%%%%%
\subsection{The example of logarithmic potentials}
%%%%%%%%%%%%%%%%%%%%%%%%%%%%%%%%%%%%%%%%%%%%
Theorem \ref{largedeviation} may directly provide some bounds where $V$ has logarithmic singularity of the type \eqref{logpotential}. As a matter of fact, one may prove coercivity  (by Lemma \ref{duality}) for any $|V|\leq \log \frac{1}{|x|}+C$,
\begin{equation}
  \int_{\T^{dN}} d\rho_N \int_{\T^{2d} \cap \{ x \ne y \} } V(x-y)\,d\mu_N^{\otimes^2}\leq \frac{1}{\kappa}\,\mathcal{H}_N(\rho_N\,\vert\;\bar\rho_N)+C_d\,\|\bar\rho\|_{L^\infty},
\end{equation}
for any $\kappa<d$, which directly implies \eqref{regularityF} for any $F$ of the form \eqref{specialF} with logarithmic singularities in $V$. In particular for $V$ as \eqref{logpotential}, $\alpha = \gamma \lambda/\kappa$ in \eqref{regularityF} and $\alpha = c/\kappa$ for the following $F_\eta$. 

We however have to be careful when checking that $I(F)=0$. The issue is then not only the singularity at $0$ of $\log |x|$ but also long distance interactions. For this reason, we truncate and denote
    \[
F_\eta(\mu)= c \,\int_{ \Pi^{2d} \cap \{x\neq y\}} \log \frac{1}{|x-y|}\,\chi_\eta(|x-y|)\,(\mu(dx)-\bar\rho(x)\,dx)\,(\mu(dy)-\bar\rho(y)\,dy),
  \]
  for $\chi_\eta(x)$ a smooth truncation function, $\chi_\eta(x)=0$ if $|x|\geq 2\eta$ while $\chi_\eta(x)=1$ if $|x|\leq \eta$. Then by employing the logarithmic Hardy-Littlewood-Sobolev inequality which has long been recognized as a key estimate for the Patlak-Keller-Segel system (see \cite{DoSe} for example),  we can show that for any $c  <d$, there exists some $\eta<1$ for which $I(F_\eta)=0$.
  As a consequence, if $V$ satisfies \eqref{logpotential}, we do not obtain a bound on $\mathcal{K}_N$ alone but instead, defining
  \[
W(x)=V(x)\,(1-\chi_\eta(x)),
\]
we obtain that there exists $C_\kappa>0$ and $\theta>0$ independent on $N$ such that
\begin{equation}
\mathcal{K}_N(G_N\,\vert\;G_{\bar\rho_N}) -\mathcal{K}_N^W(G_N\,\vert\;G_{\bar\rho_N}) \geq -\frac{\lambda}{\kappa\,\sigma}\,\mathcal{H}_N(\rho_N\,\vert\;\bar\rho_N) -\frac{C_\kappa}{N^\theta},\label{modulatedenergyPKS}
\end{equation}
as long as $\lambda<d\,\sigma$ and for any $\kappa<d$. In dimension $2$ this would correspond to imposing a total mass less than $4\,\pi\,\sigma$ in the Patlak-Keller-Segel system instead of the optimal $8\,\pi\,\sigma$.
  The issue stems from the use of general results like Theorem \ref{largedeviation} where we essentially pay twice the coercivity. Instead in our coming article \cite{BrJaWa1}, we will directly study the minimum of $\mathcal{K}_N-\mathcal{K}_N^W$ and expect to recover the optimal $\lambda<2\,d\,\sigma$ corresponding to the critical mass $8\,\pi\,\sigma$.   
%%%%%%%%%%%%%%%%%%%%%%%%%%%%%%%%%%%%%%%%%%%%%%%%
\subsection{The help from repulsive potential $V$}
%%%%%%%%%%%%%%%%%%%%%%%%%%%%%%%%%%%%%%%%%%%%%%%%%%
In the previous two subsections, the contribution $\mathcal{K}_N$ from the potential in $E_N$ was a problem and had to be controlled through the relative entropy $\mathcal{H}_N$. But if the potential is repulsive in an appropriate sense, we should expect $\mathcal{K}_N$ to be non-negative (or almost non-negative up to vanishing corrections in $N$) and even to help in bounding the right-hand side in Prop. \ref{modulatedweak} or Lemma \ref{smoothcorrection}.

This is the  strategy followed in \cite{Du,Se,Se1} for Riesz potentials $V= {C}/{|x|^\alpha}$ 
with $0<\alpha<d$. While some of the arguments in those articles rely on the explicit form of $V$ (in particular those using the extension representation in \cite{CafSil}), it is natural to ask what minimal assumptions could be required on the potential $V$, especially in the present setting where some ``defects'' in the potential may be compensated through the use of the relative entropy.  We defer to our coming article \cite{BrJaWa1} for a more thorough discussion of this issue and instead we just give simple examples of assumptions on $V$ that still allow in our setting both to bound from below $\mathcal{K}_N$ and to control the right-hand side.

The first point is the bound from below on $\mathcal{K}_N$. Of course if $\hat V\geq 0$, then one immediately has
\[
\int_{\T^{2d}} V(x-y)\,(d\mu_N-d\bar\rho)^{\otimes^2}\geq 0.
\]
But  we actually only integrate on $\T^{2d}\cap \{x\neq y\}$ removing the diagonal. This is necessary if $V(0)=+\infty$ as otherwise the self-interaction in $\mu_N\otimes\mu_N$ would immediately lead to $\mathcal{K}_N=+\infty$. Instead one may introduce truncations as in \cite{Se} and obtain

\smallskip

\begin{lemma}
Assume that $V\in L^1$, $\hat V\geq 0$, $\lim_{x\to 0} V(x)=+\infty$, $|\nabla V(x)|\leq C\,|x|^{-k}$ for some exponent $k$, and that $V$ satisfies a sort of doubling property for some constant $C>0$
  \begin{equation}
V(x)\leq C\,V(y),\qquad \mbox{if}\quad |y|\leq 2\,|x|.\label{doubling}
  \end{equation}
  Then there exists an exponent $\theta>0$ and a constant $C>0$ independent of $N$
 s.t.
  \[
\int_{\T^{2d}\cap\{x\neq y\}} V(x-y)\,(d\mu_N-d\bar\rho)^{\otimes^2}\geq -\frac{C}{N^\theta}.
  \]\label{controldiagonal}
\end{lemma}

There remains the question of the control of the right-hand side in Prop. \ref{modulatedweak}, Lemma \ref{smoothcorrection} and more generally of terms like
\[
\int_{\T^{dN}}d\rho_N\int_{\T^{2d}\cap\{x\neq y\}} \nabla V(x-y)\,(\psi(x)-\psi(y))\,(d\mu_N-d\bar\rho)^{\otimes^2},
\]
for some smooth $\psi$. Instead of using the extension representation in \cite{CafSil}, one may study this term in Fourier and with the help of the relative entropy obtain results like

\smallskip

\begin{lemma} 
  Assume that $\log\bar\rho\in W^{1,\infty}(\T^d)$ and $\psi\in W^{d+2,1}(\T^d)$  and that $V$ satisfies the assumptions of Lemma \ref{controldiagonal} with in addition 
  for any $\xi$
  \begin{equation}
|\nabla_\xi \hat V(\xi)|\leq C\,\frac{\hat V(\xi)}{1+|\xi|}+\frac{C}{1+|\xi|^{d+1}}.\qquad\label{Vhatassumpt}
\end{equation}
Then there exists a constant $\bar C>0$ depending on $C$, $\|\psi\|_{W^{d+2,1}(\T^d)}$, the constants in Lemma \ref{controldiagonal} and an exponent $\theta>0$  s.t. 
\[
\int_{\T^{dN}}d\rho_N\int_{\T^{2d}\cap\{x\neq y\}} \nabla V(x-y)\,(\psi(x)-\psi(y))\,(d\mu_N-d\bar\rho)^{\otimes^2}\leq C\,E_N\left(\frac{\rho_N}{G_N}\,\vert\;\frac{\bar\rho_N}{G_{\bar\rho_N}}\right) +\frac{C}{N^\theta}.
\]\label{repulsiverighthandside}
\end{lemma}
Lemma \ref{repulsiverighthandside} directly applies if $V={C}/{|x|^\alpha}$ where $0<\alpha <d$. It also allows to use more general potentials though again more general forms of potential will be discussed in \cite{BrJaWa1}. 
%%%%%%%%%%%%%%%%%%%%%%%%%%%%%%%%%%%%%%%%%%%%%%%%
%%%%%%%%%%%%%%%%%%%%%%%%%%%%%%%%%%%%%%%%%%%%%%%%%
\section{Towards a general result\label{general}}
%%%%%%%%%%%%%%%%%%%%%%%%%%%%%%%%%%%%%%%%%%%%%%%%%%
%%%%%%%%%%%%%%%%%%%%%%%%%%%%%%%%%%%%%%%%%%%%%%%%%%
We summarize the main ideas and bounds by sketching an example of possible general result below. Essentially we may consider even potentials $V$
\begin{equation}
V(-x)=V(x),\qquad V=V_a+V_r+V_s,\label{Vsum}
\end{equation}  
which are a sum of
\begin{itemize}
\item Any smooth enough potential $V_s$ with
  \begin{equation}
V_s\in W^{2,\infty}(\T^d);\label{smoothpotential}
  \end{equation}
\item Any possibly very singular but {\em repulsive} potential $V_r$ with 
  \begin{equation}
    \begin{split}
& \hat V_r\geq 0,\quad \lim_{x\to 0} V_r(x)=+\infty,\quad |\nabla V_r(x)|\leq C\,|x|^{-k}
    \hbox{ for some exponent } k, \\
& V_r(x)\leq C\,V_r(y),\qquad \mbox{if}\quad |y|\leq 2\,|x|,\\
&  |\nabla_\xi \hat V_r(\xi)|\leq C\,\frac{\hat V_r(\xi)}{1+|\xi|}+\frac{C}{1+|\xi|^{d +1}};
    \end{split}\label{repulsivepotential}
  \end{equation}
  \item An attractive but at most mildly singular and small part $V_a$ with for some constant $c_d$
    \begin{equation}
      \begin{split} 
        & |V_a(x)|\leq c_d+\gamma\,\log \frac{1}{|x|}\quad\mbox{for some}\ 0\leq \gamma<2d \sigma ,\quad |\nabla V_a(x)|\leq \frac{C}{|x|}, \label{attracted}\\
        & \sup_{\mu\in L^1(\T^d)} -c\,\int_{\T^{2d}} V_a(x-y)\,(d\mu-d\bar\rho)^{\otimes^2}-\sigma\,\int_{\T^d} d\mu\,\log\frac{\mu}{\bar\rho}=0, \\
        & \hskip5cm 
        \qquad\mbox{for some}\ c\ \mbox{with}\, \frac{\gamma}{2  d \sigma  } + \frac{\gamma}{2 c } < 1 . 
        \end{split}
      \end{equation}
\end{itemize}
This class contains in particular potentials with logarithmic singularities like \eqref{logpotential} provided that $\lambda<2\,\sigma$. But in general a further analysis would be required to make more explicit which potentials are allowed, while at the same time we will further improve on some of the assumptions on $V_r$ and $V_a$ in \cite{BrJaWa1}. 

\medskip

\noindent Within this class of potentials, we have the explicit propagation of chaos

\smallskip

\begin{theorem}\label{Main}
Let $\sigma>0$ and let the potential $V$ satisfies \eqref{Vsum}--\eqref{attracted}. Let $\bar\rho$
where $\inf\bar\rho>0$ solving Equation \eqref{MFlimit} with  $\bar\rho \in L^\infty\bigl(0,T;W^{2,\infty}(\T^d)\cap W^{d+3,1}(\T^d)\cap \mathcal{P}(\T^d)\bigr)$. There exists a constant $C>0$  and an exponent $\theta>0$ independent on $N$ s.t.
for  any entropy solution $\rho_N$ to Equation \eqref{liouvilleN}, we have that
\[
E_N\left(\frac{\rho_N}{G_N}\,\vert\; \frac{\bar\rho_N}{G_{\bar\rho_N}}\right)(t)\leq e^{C\,t}\, \left(E_N\left(\frac{\rho_N}{G_N}\,\vert\; \frac{\bar\rho_N}{G_{\bar\rho_N}}\right)(t=0)+\frac{C}{N^\theta}\right),
\]
which furthermore leads to an explicit bound on the relative entropy with
\[
E_N\left(\frac{\rho_N}{G_N}\,\vert\; \frac{\bar\rho_N}{G_{\bar\rho_N}}\right)(t)\geq \frac{1}{C}\,\mathcal{H}_N\left(\rho_N\,\vert\; \bar\rho_N\right)(t) -\frac{C}{N^\theta}.
\]\label{maintheorem}
\end{theorem}
All constants and exponents in Theorem \ref{maintheorem} can be made explicit with explicit dependence on the various norms of $\bar\rho$ and $V$ that are assumed to be finite.
   We finally stress that the scaled relative entropy $\mathcal{H_N}$ controls the relative entropy of all marginals with
\[
\frac{1}{k}\int_{\T^{dk}} d\rho_{N,k}\,\log \frac{\rho_{N,k}}{\bar\rho^{\otimes^k}} \leq \mathcal{H}_N\left(\rho_N\,\vert\; \bar\rho_N\right).
\]
By the Csisz\'ar-Kullback-Pinsker inequality, we hence have that
\[
\|\rho_{N,k}-\bar\rho^{\otimes^k}\|_{L^1}\leq \, \sqrt{ 2\,k}\,\left(\mathcal{H}_N\left(\rho_N\,\vert\; \bar\rho_N\right)\right)^{1/2}
\]
and Theorem \ref{maintheorem} indeed implies Theorem \ref{PKSestimate}.
% etc, etc

% The Appendices part is started with the command \appendix;
% appendix sections are then done as normal sections
% \appendix

% \section{}
% \label{}

% The Acknowledgements are an un-numbered section
%\section*{Acknowledgements}
% Acknowledgements text here
%\bibliographystyle{plain}
%\bibliography{Ref-BD}

\end{document}